\documentclass[11pt,twoside]{article}
\usepackage{amsfonts,latexsym,rawfonts,amsmath,amssymb}

\newcommand{\s}{\,\,\,\,}
\newcommand{\bthm}[2]{\vskip 8pt\bf #1\hskip 2pt\bf#2\it \hskip 8pt}
\newcommand{\ethm}{\vskip 8pt\rm}

\def\dint{\displaystyle{\int}}
\def\d#1#2{\frac{\displaystyle #1}{\displaystyle #2}}

\title{{\bf\Large{Bubbling location for $F$-harmonic maps and Inhomogeneous
Landau-Lifshitz equations}} \footnote{The second author (Y. Wang)
is supported in part by the National Key Basic Research Fund
G1999075107 and the National Science Fund for Distinguished Young
Scholars 10025104 of the People's Republic of China.}}

\author{Yuxiang Li, \s\s Youde Wang}
\date{}

\begin{document}
\maketitle

\begin{abstract}

Let $f$ be a positive smooth function on  a close Riemann surface (M,g).  The $f-energy$
of a map $u$ from $M$ to a Riemannian manifold $(N,h)$ is defined as
$$E_f(u)=\int_Mf|\nabla u|^2dV_g.$$
In this paper, we will study the blow-up properties of Palais-Smale sequences for $E_f$. We will
show that, if a Palais-Smale sequence is not compact, then it must
blows up at some critical points of $f$. As a sequence, if an
inhomogeneous Landau-Lifshitz system, i.e. a solution of
$$u_t=u\times\tau_f(u)+\tau_f(u),\s u:M\rightarrow S^2$$
blows up at time $\infty$, then the blow-up points must be the critical points of $f$.\\


\end{abstract}

\section{Introduction}
Let $(M, g)$ and $(N, h)$ be two Riemannian manifolds. A
$C^1$-smooth  map $u$ from $M$ into $N$ is called a harmonic map
if and only if $u$ is a critical point of the energy functional
$E(v)$, which is defined in local coordinates by
 $$E(v)\equiv\dint_M {\rm Trace}_g (v^*h)dV_g,$$
where
 $${\rm Trace}_g (v^*h) = g^{ij}\frac{\partial u^\alpha}{\partial x^i}\frac{\partial u^\beta} {\partial
x^j} h_{\alpha\beta}(u).$$ It is well-known that the energy functional
is conformally invariant when $\dim(M)=2$.

In this paper we would like to study a class of $C^1$-smooth maps
from a Riemann surface into a compact Riemannian manifold which
are defined as the critical points of the inhomogeneous energy
functional written as
 $$E_f(v)\equiv \dint_M {\rm Trace}_g (v^*h)fdV_g,$$
where $f$ is a smooth real function. In \cite{L} and \cite{E-L}
(see page 48, (10.20)), such maps are called $f$-harmonic from $M$
into $N$. Obviously, they are just harmonic maps if $f \equiv 1$.
Moreover, when $m=dim(M)\neq 2$, an $f$-harmonic map is nothing
but a harmonic map from $(M,f^{\frac{2}{n-2}}g)$ to $(N,h)$. In
local coordinates, the $f$-harmonic map satisfies the following
Euler-Lagrange equation
 $$f\tau(u) + \nabla f\cdot \nabla u = 0.$$
Here $\tau(u)$ is the tension field of $u$ which can be written as
 $$\tau^\alpha(u) = \Delta_g u^\alpha + g^{ij}\Gamma^\alpha_{\beta\gamma}(u)
 \frac{\partial u^\beta}{\partial x^i}\frac{\partial u^\gamma} {\partial x^j}.$$

To see the physical motivation for the $f$-harmonic maps, we
consider a smooth domain $\Omega$ in the Euclidean space
$\mathbb{R}^m$. An inhomogeneous Heisenberg spin system is given
by
 $$\partial_tu = f(u\wedge\Delta u) + \nabla f\cdot(u\wedge\nabla u),$$
where $f$ is a real-valued function defined on $\Omega$, $u(x,t)\in
S^2$, $\wedge$ denotes the cross products in $\mathbb{R}^3$ and
$\Delta$ is the Laplace operator on $\mathbb{R}^m$. Physically, the
function $f$ is called the coupling function, and is the continuum
limit of the coupling constants between the neighboring spins. It
is easy to see that if $u$ is a smooth stationary solution of the
above equation, then $u$ is just an $f$-harmonic map from $\Omega$
into $S^2$. Indeed, in this case the tension field of $u$ can be
written as $\Delta u + |\nabla u|^2u$, therefore, the right hand side
of the above equation can be expressed by $u\wedge(f\tau(u) +
\nabla f\cdot\nabla u)$, and $u$ satisfies the following equation
 $$f\tau(u) + \nabla f\cdot\nabla u = 0.$$

The above inhomogeneous Heisenberg spin system is also called
inhomogeneous Landau-Lifshitz system. Landau and Lifshitz also
suggested  considering the following dispersive system
 $$\partial_tu = u \wedge (f\tau(u) + \nabla f\cdot\nabla u) -
 u\wedge(u\wedge(f\tau(u) + \nabla f\cdot\nabla u)),$$
with an initial value condition
 $$u(0) = u_0.$$
For the well-known equation, Tang \cite{T} proved that it admits a
global weak solution which is smooth except for finitely many
points, if the domain manifold $M$ is $2$-dimensional closed , $f$
is a smooth positive function and the initial value map belongs to
$W^{1,2}(M, S^2)$ (see also \cite{St} and \cite{G-H}). The bubbles
which the weak solution blows are called as the magnetic bubbles
(\cite{Sh}). A natural question arises, {\it where do the bubbling
points of the Landau-Lifshitz equation locate}? In this paper, we
intend to answer this problem partially.

Throughout this paper, we will always assume that $f$ is  smooth
and positive. In order to answer the above question,
mathematically we need to consider the convergence and bubbling of
the sequence of $f$-harmonic maps with coupling function $f$.
Precisely, we obtain the following results.

\bthm {Theorem} {1} {\it Let D be the unit disc in $\mathbb{R}^2$.
If $u: D\setminus\{0\}\rightarrow N$ is a $W^{2,2}_{loc}$-map with
finite energy and satisfies the following equation:
$$\tau(u)=\alpha\nabla u + g$$
where $\alpha \in C^0(D)$, and $g\in L^p(D,TN)$ for some $p>2$,
then $u$ may be extended to a map $\tilde{u} \in
W^{2,p}(D,N)$.\ethm

\bthm{Theorem}{2} Let $(M,g)$ be a closed Riemann surface, and $N$
 a compact submanifold of $\mathbb{R}^K$. Let $f$ be a smooth
positive function on $M$. Assume that $u_k\in W^{2,2}(M,N)$ is a
sequence which satisfies
 $$f\tau(u_k) + \nabla f\cdot\nabla u_k = \alpha_k$$
and
 $$\int_M|\nabla u_k|^2fdV_g\leq C,$$
where $\alpha_k\in L^2(u_k^{-1}(TN))$ and satisfies
$$||\alpha_k||_{L^2}\rightarrow 0,\,\,\, as\,\, k\rightarrow+\infty$$
If $p$ is a blow-up point of the sequence, i.e.
$$\lim_{r\rightarrow 0}\liminf_{k\rightarrow+\infty}
\dint_{B_r(p)}|\nabla u_k|^2fdV_g>0,$$ then, $p$ must be a
critical point of $f$. \ethm

Applying the above theorem to the inhomogeneous Landau-Lifshitz
equation, we can partially answer  the question mentioned before.
Concretely, we come to the conclusions

\bthm{Theorem}{3} Let $(M,g)$ be a closed Riemann surface, and
$S^2$  a unit sphere with standard metric. Suppose that the
coupling function $f$ is smooth and positive  on $M$ and $u \in
L^2((0, \infty); W^{2,2}(M, S^2))$ is the unique weak solution for
the initial value problem of the inhomogeneous Landau-Lifshitz
equation with initial map $u_0\in W^{1,2}(M, S^2)$. If $u(t)
\equiv u(\cdot, t)$ blows up at time infinity, then, the blow-up
points must be the critical points of the coupling function $f$.
\ethm

\section{Removable singularity}

It is well-known that the removable singularity theorem of
Sacks-Uhlenbeck says that a harmonic map from
$D\setminus\{0\}\rightarrow N$ with finite energy can be extend to
0 smoothly. The main aim of this section is to generalize
Sacks-Uhlenbeck's theorem to the present case, i.e., to prove
Theorem 1. The method adopted here is due to Sacks-Uhlenbeck
essentially. One still sees that the Hopf differential is the key
in the proof. However, in our case, the Hopf differential is no
longer holomorphic, thus the proof should be a little more
delicate than theirs.

Let us first recall the $\epsilon$-regularity discovered by Sacks
and Uhlenbeck.

\bthm{Lemma}{2.1} Suppose that $u\in W^{2,2}(D,N)$ satisfying
$$\tau(u)=g\in L^2(D,TN).$$
Then there exits $\epsilon > 0$ such that if $\int_{D}|\nabla
u|^2\leq \epsilon$ we have
$$||u-\bar{u}||_{W^{2,2}(D_{\frac{1}{2}})}
\leq C(||\nabla u||_{L^2(D)} +||g||_{L^2(D)}).$$ Here $\bar{u}$ is
the mean value of $u$ over the unit disc and $D_{\frac{1}{2}}$ is
a disc with radius $\frac{1}{2}$ and centered at the origin.\ethm

{\bf Proof:} cf \cite{S-U}, or \cite{D},or \cite{D-T}.

\vskip2mm

Using the standard elliptic estimate, we have

\bthm{Corollary}{2.2} Suppose that $u\in W^{2,2}(D,N)$ satisfies
 $$\tau(u)=\alpha\nabla u + g,\eqno (2.1)$$
where $\alpha(x)\in C^0(D)$ and $g\in L^p(D,TN)$ for some $p>2$.
Then, there exists $\epsilon>0$ such that whenever
$\int_{D}|\nabla u|^2\leq \epsilon$ we have
$$|\nabla u|(0)\leq C(||\alpha||_{C^0(D)},p)(||\nabla u||_{L^2(D)} +
||g||_{L^p(D)}).$$ \ethm

In this section, we always assume $u$ to be a map from
$D\setminus\{0\}$ to N which belongs to
$W^{2,2}_{loc}(D\setminus\{0\}, N)$ and satisfies the equation
(2.1). In order to prove Theorem 1, we need to prove the following
lemmas. First, we have

\bthm {Lemma} {2.3} There exists $\epsilon>0$ such that if
$\int_{D}|\nabla u|^2dx < \epsilon$, then there holds true
$$|x||\nabla u|(x)\leq C(||\nabla u||_{L^2(D_{2|x|})}
+|x|^{2-\frac{2}{p}}\times||g||_{L^p(D_{2|x|})})\s \mbox{for} \s
\forall x\in D_{\frac{1}{2}},$$ where $C$ is a positive constant
which depends only on $\epsilon$.\ethm

{\bf Proof:} Fix an $x_0\in D_{\frac{1}{2}}$, we define
$\tilde{u}=u(x|x_0|+x_0)$. Then we have
$$\tau(\tilde{u})=|x_0|^2g+|x_0|\nabla\tilde{u}.$$
Notice that $|\nabla\tilde{u}|(0)=|\nabla u||x_0|$, then we can
get this Lemma from Corollary 2.2.

\vskip1.5cm

Now, let
 $$\Psi= \langle u_x,u_x\rangle - \langle u_y,u_y\rangle - 2i\langle u_x,u_y\rangle
 =4\langle\frac{\partial u}{\partial z}, \frac{\partial u}{\partial
 z}\rangle,$$
where $z=x + iy$. It is easy to see that
 $$\partial_{\bar{z}}\Psi=8\langle \Delta u,\frac{\partial u}{\partial
 z}\rangle = 8\langle A(u)(du,du)+\alpha\nabla u + g,\frac{\partial u}{\partial z}\rangle
 = 8 \langle \alpha \nabla u+g, \frac{\partial u}{\partial z}\rangle.\eqno (2.2)$$ We
need to prove a Stokes type equality for the 1-form $z\Psi$.

\bthm {Lemma} {2.4} There holds true that
$$\int_{|z|=r}z \Psi dz = \int_{D_r}z\partial_{\bar{z}}
\Psi d\bar{z}\wedge dz.$$ \ethm

{\bf Proof:} As
 $$ d(z \Psi dz) = \partial_{\bar{z}}(z\Psi dz) = z\partial_{\bar{z}}\Psi
 d \bar{z}\wedge dz,$$
by applying the Stokes formula, for any $r_0<r$, we have
$$\int_{|z|=r\setminus|z| = r_0}z\Psi dz=\int_{D_r\setminus D_{r_0}}
z\partial_{\bar{z}}\Psi d\bar{z}\wedge dz.$$ By (2.2),
 $$\int_{D_{r_0}}|z\partial_{\bar{z}}\Psi d\bar{z}\wedge
 dz| \leq C r_0\int_{D_{r_0}}(|\alpha\nabla u|^2+|g|^2)dx \rightarrow
 0$$
as $r_0\rightarrow 0$. Therefore, to complete the proof of the
Lemma, we only need to prove
 $$\int_{|z|=r_0}z\Psi dz=\sqrt{-1}\int_0^{2\pi}z^2\Psi
 d\theta|_{|z|=r_0}\rightarrow 0.$$
However the last equality follows from Lemma 2.3.

\bthm {Lemma}{2.5} There holds
$$\int_{D_r}\langle u_r,u_r\rangle -\int_{D_r}\langle u_\theta,u_\theta\rangle
=O(r).$$\ethm

{\bf Proof:} By a direct computation, we have
$$Re (z^2\Psi)=-|u_\theta|^2+|z|^2|u_r|^2,$$
then
$$\begin{array}{ll}
 |Re\dint_0^{2\pi}z^2\Psi d\theta|&=|Im\dint_{|z|=r}z\Psi
 dz|\\[1.5ex]
&=|Im\dint_{D_r}z\partial_{\bar{z}}\Psi dz\wedge d\bar{z}|\\[1.5ex]
 &\leq \dint_{D_r}|z\partial_{\bar{z}}\Psi dz\wedge d\bar{z}|\\[1.5ex]
    &\leq r\dint_{D_r}|\alpha\nabla u+g||\nabla u|dx\\[1.5ex]
    &\leq C r\dint_{D_r}(|g|^2+|\alpha\nabla u|^2)dx.
  \end{array}$$
i.e.
$$\int_0^{2\pi}|u_r(r,\theta)|^2r^2d\theta-
\int_0^{2\pi}|u_\theta(r,\theta)|^2d\theta=O(r).$$ Therefore
$$\begin{array}{ll}
    \dint_{D_r}(\langle u_r,u_r\rangle - \langle u_\theta,u_\theta \rangle) &=
          \dint_0^r\int_0^{2\pi}(|u_r|^2-\d{1}
           {r^2}|u_\theta|^2)rd\theta dr\\[1.5ex]
      &=\dint_0^r\d{1}{r}O(r)dr = \dint_0^rO(1)dr\\[1.5ex]
      &=O(r).
  \end{array}$$

{\bf Proof of Theorem 1:} As in \cite{S-U}, we approximate u by
the function $q$ which is harmonic on every domain,
$$D_m(r_0)=\{z:2^{-m-1}r_0< |z|< 2^{-m}r_0\},$$
and equals to
$$\frac{1}{2\pi}\int_0^{2\pi}u(2^{-m}r_0,\theta)
d\theta$$ and
$$\frac{1}{2\pi}\int_0^{2\pi}u(2^{-m-1}r_0,\theta)
d\theta$$ respectively on the boundaries $\{z: |z|=2^{-m}r_0\}$
and $\{z: |z|=2^{-m-1}r_0\}$. Then $q$ is piecewise linear in
$\log{r}$ and depends only on the radial coordinate. Now, for
$2^{-m-1}r_0 \leq r \leq 2^{-m}r_0$,
$$\begin{array}{ll}
   |q(r)-u(r,\theta)|&\leq |q(r)-q(2^{-m}r_0)|+|q(2^{-m}r_0)-u(r,\theta)|\\[1.5ex]
      &\leq |q(2^{-m-1}r_0)-q(2^{-m}r_0)|+
         \d{1}{2\pi}\dint_0^{2\pi}|u(2^{-m}r_0,\theta^\prime)-
         u(r,\theta)|d\theta^\prime\\[1.5ex]
      &\leq C\sup_{2^{-m-1}r_0\leq|r|\leq 2^{-m}r_0} r|\nabla u|\\[1.5ex]
          &\leq C(||\nabla u||_{L^2(D_{2r})}+
          r^{2-\frac{2}{p}}||g||_{L^p(D_{2r})}).
  \end{array}$$

Now, we estimate the difference between $q$ and $u$:
$$\begin{array}{ll}
    \dint_{D_{r_0}}|\nabla (u-q)|^2&=\sum_{m=0}^\infty
      r\dint_0^{2\pi}(u(r,\theta) - q(r))(u_r(r,\theta)-q'(r))d\theta
       |_{2^{-m-1}r_0}^{2^{-m}r_0}\\[1.5ex]
       &\;\;-\dint_{D_r}
      (q-u)\Delta(q-u)dx.
  \end{array}\eqno (2.3)$$
Since $q'(r)=constant\times\frac{1}{r}$ on $D_m(r_0)$,
$$\int_0^{2\pi}(u(r,\theta) - q(r))q'(r)d\theta=0,\;\;\forall
r=2^{-m}r_0,$$ Hence,
$$\begin{array}{ll}
   &\sum_{m=0}^\infty r\dint_0^{2\pi}(u(r,\theta) - q(r))
      (u_r(r,\theta)-q'(r))d\theta|_{2^{-m-1}r_0}^{2^{-m}r_0}\\[1.5ex]
    =& r_0\dint_0^{2\pi}(u(r_0,\theta) - q(r_0))u_r(r_0,
         \theta)d\theta\\[1.5ex]
     &-\lim\limits_{m\rightarrow+\infty}2^{-m}r_0\dint_0^{2\pi}
      (u(2^{-m}r_0,\theta) - q(2^{-m}r_0))u_r(2^{-m}r_0,\theta))d\theta.
 \end{array}$$
By Lemma 2.3, we have
$$\begin{array}{ll}
       &2^{-m}r_0\dint_0^{2\pi}(u(2^{-m}r_0,\theta) - q(2^{-m}r_0))u_r(2^{-m}r_0,
       \theta)d\theta\\[1.5ex]
   \leq & ||u(2^{-m}r_0,\theta)-q(2^{-m}r_0)||_{L^\infty}\sup_{r=2^{-m}r_0} r
   |\nabla u(r,\theta)|\\[1.5ex]
   \rightarrow & 0
  \end{array}$$
as $m\rightarrow+\infty$.

Moreover, we have
$$\begin{array}{ll}
    r_0\dint_0^{2\pi}(u(r_0,\theta)-q(r_0))u_r(r_0,\theta)d\theta
      &\leq r_0\left(\dint_0^{2\pi}(u(r_0,\theta)-q(r_0))^2d\theta
        \dint_0^{2\pi}|u_r(r_0,\theta)|^2d\theta\right)^\frac{1}{2}\\[1.5ex]
      &\leq\left(\dint_0^{2\pi}|u_\theta(r_0,\theta)|^2d\theta\right)^\frac{1}{2}
        \left(\dint_0^{2\pi}r_0^2|u_r(r_0,\theta)|^2d\theta\right)^\frac{1}{2}\\[1.5ex]
      &\leq
\displaystyle\frac{1}{2}\dint_0^{2\pi}(|u_\theta(r_0,\theta)|^2+|u_r(r_0,\theta)|^2r_0^2)d\theta\\[1.5ex]
      &=\displaystyle{\frac{r_0}{2}}\dint_0^{2\pi}|\nabla u(r_0,\theta)|^2r_0d\theta.
 \end{array}\eqno (2.4)$$
and
$$\begin{array}{ll}
   \int_{D_{r_0}}|(q-u)(\Delta(q-u))|&=\dint_{D_{r_0}}|q-u|\times
       |A(u)(du,du)-\alpha\nabla u - g|dx\\[1.5ex]
     &\leq ||q-u||_{L^\infty(D_{r_0})}(
         ||A||_{L^\infty}\dint_{D_{r_0}}|\nabla u|^2dx + \sqrt{\pi}r_0||\alpha\nabla u +
g||_{L^2(D_{r_0})})
  \end{array}\eqno (2.5)$$
Obviously, for any $1>\delta>0$, we can always pick up $r_0$ which
is small enough such that
$$\int_{D_{r_0}}|(q-u)(\Delta(q-u))| \leq \delta(\int_{D_{r_0}}| \nabla u|^2dx+r_0).$$
Applying Lemma 2.5, we get
$$\begin{array}{lll}
   \dint_{D_{r_0}}|\nabla (u-q)|^2dx
     &\geq&\dint_{D_{r_0}}\langle u_\theta,u_\theta\rangle
       dx\\[1.5ex]
     &=&\d{1}{2}\dint_{D_{r_0}}(\langle u_\theta,u_\theta\rangle + \langle u_r,u_r\rangle)dx \\
     & &+ \d{1}{2}\dint_{D_{r_0}}(\langle u_\theta,u_\theta\rangle-
        \langle u_r,u_r\rangle)dx\\[1.5ex]
     &=&\d{1}{2}\dint_{D_{r_0}}|\nabla u|^2dx+O(r_0).
  \end{array}\eqno (2.6)$$
Then, from (2.3), (2.4), (2.5) and (2.6) we can derive
 $$\lambda\int_{D_{r_0}}|\nabla u|^2\leq r_0\int_0^{2\pi}|\nabla u(r_0
 ,\theta)|^2r_0d\theta + Cr_0 ,$$
where $\lambda$ is a positive constant which is smaller than 1.

Set
 $$f(r)=\int_{D_r}|\nabla u|^2dx,$$
then we have
 $$\lambda f(r)<rf'(r)+Cr .$$
Hence,
 $$(\frac{f}{r^\lambda})'\geq -Cr^{-\lambda} .$$
By integrating the above differential inequality over the interval
$[r, 1]$, we obtain
 $$f(r)\leq Cr^\lambda\int_r^1 s^{-\lambda}ds + f(1)r^\lambda \leq Cr^{\lambda}.$$
By applying Lemma 2.3, it follows from the above inequality that
$$|\nabla u|(x)\leq |x|^{\lambda-1}.$$
Thus, we can complete the proof of the theorem by standard
elliptic estimate theory.

\section{A variational formula}

For the inhomogeneous functional $E_f(\cdot)$ defined on
$W^{1,2}(M,N)$,  we can see easily that the first variational
formula at point $u\in W^{2,2}(M,N)$ can be written as
 $$dE_f(\xi)=\int_M<f\tau(u)+\nabla u\nabla f,\xi>dV_g,$$
for any $\xi\in T_uW^{1,2}(M,N)$. Here, we need to derive another
formula of $E_f(\cdot)$ with respect to the variation of the
domain manifold. The following calculation is essentially due to
Price (\cite{P}).

Take an $1$-parameter family of transformations $\{\phi_s\}$ of M,
which is generated by the vector field $X$, we have
$$\begin{array}{ll}
   E_f(u\circ\phi_s)&=\displaystyle\frac{1}{2}\dint_M|\nabla (u\circ\phi_s)|^2f(x)dV_g\\[1.5ex]
      &=\displaystyle\frac{1}{2}\dint_M\sum_\alpha|d(u\circ\phi_s)(e_\alpha)|^2f(x)dV_g(x)\\[1.5ex]
      &=\displaystyle\frac{1}{2}\dint_M\sum_\alpha|du({\phi_{s}}_*(e_\alpha))
       |^2(\phi_s(x))f(x)dV_g(x)\\[1.5ex]
     &=\displaystyle\frac{1}{2}\dint_M\sum_\alpha|du({\phi_{s}}_*(e_\alpha))|^2(x)
      f(\phi_{-s})Jac(\phi_s^{-1})dV_g,
  \end{array}$$
where $\{e_\alpha\}$ is a local orthonormal basis of $TM$. Noting
 $$\frac{d}{ds}Jac(\phi_s^{-1})dV_g|_{s=0}=-div(X)dV_g,\s\s \frac{d}{ds}
 f(\phi_{-s})=-df(X),$$
we have
 $$\begin{array}{ll}
 \displaystyle\frac{d}{ds}E_f(u\circ\phi_s)|_{s=0}= &-\displaystyle\frac{1}{2}\dint_M
 |\nabla u|^2fdiv(X)dV_g-\displaystyle\frac{1}{2}\dint_Mdf(X)|\nabla u|^2dV_g\\
 & + \sum_{\alpha}\dint_M\langle du(\nabla_{e_\alpha}X), du(e_\alpha)\rangle
fdV_g.\end{array}$$ So, we have proved the formula
 $$\begin{array}{ll}
 dE_f(u)(u_*(X)) = &-\displaystyle\frac{1}{2}\dint_M|\nabla u|^2fdiv(X)dV_g
 - \displaystyle\frac{1}{2}\dint_Mdf(X)|\nabla u|^2dV_g\\
 &+\sum_{\alpha}\dint_M\langle du(\nabla_{e_\alpha} X),du(e_\alpha)\rangle
fdV_g.\end{array}$$

\section{The proof of Theorems}

The task of this section is to prove our theorem 2 and 3. In fact,
what Theorem 2 concerns is just the blow-up analysis for a
so-called Palais-Smale sequence of $E_f(u)$. We will focus on what
occur if the sequence is not compact in the Sobolev space
$W^{1,2}(M, N)$.

{\bf Proof of Theorem 2:} By the assumptions stated in Theorem 2,
$\{u_k\}\subset W^{2,2}(M,N)$ is a Palais-Smale sequence of maps
from $M$ into $N$. Then, it satisfies
 $$f\tau(u_k)+\nabla f\nabla u_k=\alpha_k \eqno (4.1)$$
 $$E_f(u_k)\leq C,\eqno (4.2)$$
where $\alpha_k\in u_k^{-1}(TN)$ and satisfies
 $$||\alpha_k||_{L^2}\rightarrow 0.\eqno(4.3)$$

First, we note that in a local complex coordinates (4.1) can be
written as
 $$f\tau_0(u_k) + \nabla_0 f\cdot\nabla_0 u_k = |\beta|\alpha_k$$
where $\tau_0$ and $\nabla_0$ are the operators defined on
$\mathbb{R}^2$ with standard Euclidean metric, since $\tau$ is
conformally invariant operator. Without  loss of generality, we
may assume $g=dx^2+dy^2$ on a complex coordinate system.

Set
 $${\cal S}\equiv \{x: \lim_{r\rightarrow 0}\liminf_{k\rightarrow +\infty}
\int_{D_r(x)}
 |\nabla u_k|^2dV_g >0\}.$$
Usually, we say $x$ is a bubbling point for the sequence $\{u_k\}$
if and only if $x\in {\cal S}$. It is easy to see that
 $\cal S$ contains only finite points. By the Lemma 2.1, for any
$x_0\in \cal{S}$, we have
$$\liminf_{k\rightarrow +\infty} \int_{D_r(x)}
 |\nabla u_k|^2dV_g >\epsilon,\s \mbox{for any} \s r>0.$$

By the weak compactness of $W^{1,2}(M,N)$ we know that there
exists a subsequence of $\{u_k\}$, still denoted by $\{u_k\}$, and
$u \in W^{1,2}(M,N)$ with
$$E_f(u)< + \infty,$$
such that $\{u_k\}$ converges weakly to $u$ in $W^{1,2}(M,N)$,
which is an $f-$harmonic map. Moreover, Theorem 1, Lemma 2.1,
Corollary 2.2 and elliptic estimate theory tell us that $u\in
C^\infty(M,N)$ and
$$u_k\rightarrow u$$
in $W^{1,q}(\Omega,N)$ for any $\Omega\subset\subset M\setminus
{\cal S}$ and $q>1$.

Thus, to prove Theorem 2 we need only to prove
$${\cal S}\subset \{\mbox{the critical points of}\s f\}.$$
Now, pick up a point $p\in S$. As we have pointed out, we may
assume $g=dx^2+dy^2$ in a complex coordinate chart $\mathcal{N}$
around $p$. Without the loss of generality, we may assume that
$p=(0,0)$, $Q=[-1,1]\times[-1,1]\subset\mathcal{N}$ and $Q\cap
S=\{p\}$. If $p$ is not the critical point of $f$, then, without
loss of generality, we may suppose that
 $$df(0)=\lambda dx,$$
where $\lambda$ is a positive constant. Thus, in the neighborhood
of $p$, $df(x)=\lambda dx+O(r)$ where $r^2=x^2+y^2$.

We need to choose two functions to cut off the vector field
$\frac{\partial}{\partial x}$  respectively in $x$ and $y$
directions. First, we take a cut off function $\sigma\in
C^\infty(\mathbb{R})$ which is 1 on $[-\delta,\delta]$, and 0 on
$[-2\delta,2\delta]^c$, where
$$\delta=\frac{\lambda\epsilon}{16||\nabla u||_{C^0(Q)}^2||f||_{C^0}}.$$
Then, we define the other function as following
$$\eta(t)=\left\{\begin{array}{ll}
                1&\mbox{if}\s |t|\leq b'\\
                (b-t)/(b-b')&\mbox{if}\s b'\leq t\leq b\\
                (b+t)/(b-b')&\mbox{if}\s -b\leq t\leq -b'\\
                0&\mbox{if}\s t>b\;\;or\;\;t<-b.
          \end{array}\right.$$
Here $b$ and $b^\prime$ are chosen to satisfy $0<a<b^\prime <b
<2a<1$ where $a$ is a constant such that
$$\int_{[-2a,2a]\times[-1,1]}|\nabla u|^2 < \frac{\lambda\epsilon}
{8||\sigma'||_{C^0}||f||_{C^0}}. \eqno(4.4)$$ Set
 $$X=\eta(x)\sigma(y)\frac{\partial}{\partial x}.$$
By a direct computation, we have
$$div(X)=\eta'(x)\sigma(y),$$
 $$\sum_\alpha\langle du_k(\nabla_{e_\alpha} X),du_k(e_\alpha)\rangle = \eta'(x)\sigma(y)
 \left|\frac{\partial u_k}{\partial x}\right|^2 + \eta(x)\sigma'(y)\langle
 \frac{\partial u_k}{\partial x}, \frac{\partial u_k}{\partial
 y}\rangle.$$
By the formula derived in section 3, we have
$$\begin{array}{ll}
  &-\dint_Q\eta'(x)\sigma(y)(\left|\d{\partial u_k}{\partial x}\right|^2
  -\left|\d{\partial u_k}{\partial y}\right|^2)fdxdy
  -2\dint_Q\eta(x)\sigma^\prime(y)\langle \d{\partial u_k}{\partial x},
  \d{\partial u_k}{\partial y}\rangle fdxdy\\[1.5ex]
  = & \dint_Q(\lambda+O(r))\eta(x)\sigma(y)|\nabla u_k|^2dxdy+
  2\dint_Q\langle\alpha_k,{u_{k}}_*(X)\rangle dxdy.
  \end{array}$$
Note that
 $$\mbox{supp}(\eta'(x)\sigma(y))\cup
\mbox{supp}(\eta(x)\sigma'(y))\subset
 Q\setminus (-a,a)\times(-\delta,\delta),$$
we can replace $Q$ in the left side of the above equality with
$Q\setminus (-a,a)\times(-\delta,\delta)$.

For arbitrarily fixed $a$ and $\delta$, $u_k$ is bounded in
$W^{2,2}(Q\setminus(-\frac{a}{2},\frac{a}{2})\times(-\frac{\delta}{2},
\frac{\delta}{2}))$. So, by taking a subsequence, we have $\nabla
u_k\rightarrow \nabla u$ in $L^2(Q\setminus(-a,
a)\times(\delta,\delta))$. Therefore
$$-\int_Q\eta'(x)\sigma(y)(\left|\frac{\partial u_k}{\partial x}\right|^2
  -\left|\frac{\partial u_k}{\partial y}\right|^2)fdxdy
  \rightarrow -\int_Q\eta'(x)\sigma(y)(\left|\frac{\partial u}{\partial x}\right|^2
  -\left|\frac{\partial u}{\partial y}\right|^2)fdxdy,\eqno (4.5)$$
and
 $$\begin{array}{ll}
  -2\dint_Q\eta(x)\sigma'(y)\langle\frac{\partial u_k}{\partial x},
  \frac{\partial u_k}{\partial y}\rangle fdxdy
  &\rightarrow
  -2\dint_{[-2a,2a]\times[-1,1]}\eta(x)\sigma'(y)\langle\frac{\partial u}{\partial x},
  \frac{\partial u}{\partial y}\rangle fdxdy\\[1.5ex]
  &\leq 2||\sigma^\prime||_{C^0}||f||_{C^0}\dint_{[-2a,2a]\times[-1,1]}|\nabla
u|^2dxdy\\[1.5ex]
  &<\d{\lambda\epsilon}{4},
  \end{array}\eqno (4.6)$$
where we have used (4.4) in the last inequality. Moreover, once
$\delta$ and $a$ are chosen, then
$$\int_{[-a,a]\times[-\delta,\delta]}|\nabla u_k|^2dxdy\geq\epsilon$$
when $k$ is sufficiently large. Hence
$$\int_Q(\lambda+O(r))\eta(x)\sigma(y)|\nabla
u_k|^2dxdy+2\int_Q\langle\alpha_k,{u_k}_*(X)\rangle dxdy
 >\frac{1}{2}\lambda\epsilon.\eqno (4.7)$$
In view of (4.5), (4.6) and (4.7), we have
 $$\int_Q\eta'(x)\sigma(y)(-\left|\frac{\partial u}{\partial x}\right|^2
 +\left|\frac{\partial u}{\partial y}\right|^2)fdxdy
 > \frac{1}{4}\lambda\epsilon.$$
Letting $b'\rightarrow b$, we get
 $$\int_{|x|=b}\sigma(y)(-\left|\frac{\partial u}{\partial x}\right|^2
 +\left|\frac{\partial u}{\partial y}\right|^2)fdy > \frac{1}{4}\lambda\epsilon.$$
Recall that $\mbox{supp} \sigma \subset [-2\delta,2\delta]$, then
 $$\delta> \frac{\lambda\epsilon}{16||\nabla u||_{C^0}^2||f||_{C^0}}$$
which contradicts the definition of $\delta$. It means that
$\lambda$ must be zero. Therefore $p$ is a critical point of $f$.
Thus we complete the proof of Theorem 2.

\vskip3mm Now, we return back to our problem on the location of
the bubbling points of the weak solutions to the inhomogeneous
Landau-Lifshitz equations. Consider the following initial value
problem:
$$\left\{\begin{array}{l}
  \partial_tu = u\wedge(f(x)\Delta u + \nabla f\cdot\nabla u) -
u\wedge(u\wedge(f(x)\Delta u + \nabla f\cdot\nabla u)),\\
  u(0) = u_0(x)\in W^{1,2}(M, S^2).
  \end{array}\right.$$
Noting $|u|^2 \equiv 1$ and the following identity
$$u\wedge(u\wedge(f(x)\Delta u + \nabla f\cdot\nabla u)) = (u\cdot(f(x)\Delta u +
\nabla f\cdot\nabla u))u - (u\cdot u)(f(x)\Delta u + \nabla
f\cdot\nabla u),$$ we can see easily that the above equations may
be rewritten by
$$\left\{\begin{array}{l}
  \partial_tu = f(x)\tau(u) + \nabla f\cdot\nabla u + u\wedge(f(x)\tau(u) + \nabla
f\cdot\nabla u),\\
  u(0) = u_0(x)\in W^{1,2}(M, S^2).
  \end{array}\right.\eqno (4.8)$$
Here $\tau(u)=\Delta u + |\nabla u|^2u$ is the tension field of
map $u: M \rightarrow S^2$.

Tang has ever employed Struwe's method to study the existence and
uniqueness of the above equation. We outline the argument in
\cite{T} as follows

1. $\exists T>0$, s.t.  the solution (4.8) solvable on
$M\times[0,T)$.

2. $u(t)$ blow-up at finitely many points.

3. $u(t)$ converges to a $u(T)\in W^{1,2}(M, N)$ weakly, and on
any sub-domain which does not contain any bubbling points, $u(t)$
strongly converges to $u(T)$ locally.

Then, we construct a new flow which stems from $u_{T}$. Then, by
the same argument as in \cite{St} we know that There exists
$T_1>0$ such that the new flow exists on the interval $[T,T_1)$
and blows up at $T_1$. At each bubbling point, $u(t)$ blows one or
more bubbles, i.e one or more non-constant harmonic maps. It is
well-known that $u(t)$ must lose energy at every bubbling point.
Hence, we always have a $\tilde{T}$ such that
$$\left\{\begin{array}{l}
  \partial_tu = f(x)\tau(u) + \nabla f\cdot\nabla u + u\wedge(f(x)\tau(u)
  + \nabla f\cdot\nabla u),\\
  u(0) = u(\tilde{T})\in W^{1,2}(M, S^2).
  \end{array}\right.$$
is solvable on $[0,\infty)$. The results in \cite{T} can be
summarized in the following lemma.

\bthm{Lemma}{4.1} Let $(M, g)$ be a closed Riemann surface and $f$
be a smooth positive function on $M$. For any $u_0\in W^{1,2}(M,
S^2)$ there exists a distribution solution $u: M \times
\mathbb{R}^ + \rightarrow S^2$ of the above equation which is
smooth on $M\times\mathbb{R}^+$ away from at most finitely many
points $(x_k, t_k)$, $1\le k\le K_0$, $0< t_k \le \infty$, which
satisfies the energy inequality $E_f(u(s)) \le E_f(u(t))$ for all
$0\le s\le t$, and which assumes its initial data continuously in
$W^{1,2}(M, S^2)$. The solution is unique in this class.\ethm

It is easy to see the following identity holds for the solution to
(4.8) and any $0<t_1<t_2\le \infty$
$$E_f(u(t_1))-E_f(u(t_2))=-\int_{t_1}^{t_2}||\partial_tu||_{L^2}^2dt.$$
Hence, it follows
$$\int_{0}^{+\infty}||\partial_tu||^2_{L^2}<+\infty.$$
This implies that there exists a sequence $\partial_tu(x, t_i)$ such
that
 $$||\partial_t u(x, t_i)||_{L^2} \rightarrow 0.$$
So, $\{u(t_i)\}$ is a Palais-Smale sequence of $E_f(u)$.
Therefore, if $u(t)$ does blow up at infinity, then, by applying
Theorem 2 we obtain the conclusion of Theorem 3.

\vskip2mm

As another example, we may  also consider the gradient flow of the
function $E_f$, i.e. a solution of
$$\left\{\begin{array}{l}
     u_t = f(x)\tau(u) + \nabla f\cdot\nabla u,\\
     u(0) = u_0(x)\in W^{1,2}(M,N).
  \end{array}\right.$$
If $u(t)$ does blow up at infinity, we also have  results similar
to Theorem 3.

\vspace{1.5cm}

Yuxiang  Li

{ICTP, Mathematics Section, Strada Costiera 11, I-34014 Trieste,
Italy}

{\it E-mail address:} liy@ictp.it\\

Youde Wang

{Institute of Mathematics, AMSS, CAS}

{\it E-mail address:} wyd@math.ac.cn}

\end{document}